\documentclass[journal]{IEEEtran}
\usepackage{times}
\usepackage{epsfig}
\usepackage{epstopdf}
\usepackage{graphicx}
\usepackage[fleqn]{amsmath}
\usepackage{amssymb}
\usepackage{subfigure}
\usepackage{booktabs}
\usepackage[lined,ruled,linesnumbered,commentsnumbered]{algorithm2e}
\usepackage{amsthm}
\usepackage{url}
\newtheorem{theorem}{Theorem}[section]
\newtheorem{lemma}[theorem]{Lemma}

\newtheorem{definition}[theorem]{Definition}

\begin{document}

\title{Low-Rank  Structure Learning via Log-Sum Heuristic Recovery
\thanks{Y. Deng, Q. Dai and
Z. Zhang are with the Department of Automation, Tsinghua University, Beijing, China (email: daiqh@tsinghua.edu.cn).}
\thanks{R. Liu is with School of Mathematical Sciences, Dalian University of Technology,
Dalian, China. (email:  rsliu0705@gmail.com).}
\thanks{ S. Hu is with College of Computer Science, Hangzhou Dianzi University, Hangzhou, China (email: sqhu@hdu.edu.cn).}
}

\author{Yue Deng,\ \ Qionghai Dai, {\it Senior Member, IEEE},\ \ Risheng Liu,\ \
 \\Zengke Zhang, \ \ Sanqing Hu, {\it Senior Member, IEEE}\\
}
\date{}

\maketitle

\pagenumbering{arabic}

\begin{abstract}
Recovering intrinsic data structure from corrupted observations plays an important role in various tasks in the communities of machine learning and signal processing. In this paper, we propose a novel model, named log-sum heuristic recovery (LHR), to learn the essential low-rank structure from corrupted data. Different from traditional approaches, which directly utilize $\ell_1$ norm to measure the  sparseness, LHR introduces a more reasonable log-sum measurement to enhance the sparsity in both  the intrinsic low-rank structure and in the sparse corruptions. Although the proposed LHR optimization is no longer convex, it still can be effectively solved by a majorization-minimization (MM) type algorithm, with which the non-convex objective function is iteratively replaced by its convex surrogate and  LHR finally falls into the general framework of reweighed approaches.
We prove that the MM-type algorithm can converge to a stationary point after successive iteration. We test the performance of our proposed model by applying it to solve two typical problems: robust principal component analysis (RPCA) and low-rank representation (LRR).
 For RPCA, we compare LHR with the benchmark Principal Component Pursuit (PCP) method from both the perspectives of simulations and practical applications. For LRR, we apply LHR to compute the low-rank representation matrix for motion segmentation and stock clustering. Experimental results on low rank structure learning demonstrate that the proposed Log-sum based model performs much better than the $\ell_1$-based method on for data with higher rank and with denser corruptions.
\end{abstract}

\section{Introduction}
Learning the intrinsic data structures via matrix analysis \cite{li2010deterministic}\cite{yuan2009binary} has received wide attention in many fields, e.g., neural network\cite{TNN_rank}, learning system\cite{Matrix-classification}\cite{LRR}, control theory \cite{Automatic_rank}, computer vision \cite{Robust-PCA}\cite{MC_fusion} and pattern recognition
\cite{MC_PR}\cite{Deng_CVIU}. There are quite a number
of efficient mathematical tools for rank analysis, e.g., Principal
Component Analysis (PCA) and Singular Value Decomposition (SVD).
However, these typical approaches could only handle some preliminary
and simple problems. With the recent progresses of compressive
sensing \cite{Donohol-CS},  a new concept on nuclear norm
optimization has emerged into the field of rank minimization
\cite{nuclear-norm} and has led to a number of interesting applications, e.g. low rank structure learning (LRSL) from corruptions.

LRSL is a general model for many practical problems in the communities of machine learning and signal processing, which considers learning a data of
the low rank structure from sparse errors \cite{Candes:1229741}
\cite{Chandrasekaran:1183674}\cite{RMD}\cite{ICML-LRR}. Such
problem can be formulated as: $\mathbf{P} = f(\mathbf{A}) +
g(\mathbf{E})$, where $\mathbf{P}$ is the corrupted matrix
observed in practical world; $\mathbf{A}$ and $\mathbf{E}$ are \emph{low-rank} matrix
 and   \emph{sparse}
corruption, respectively and the functions $f(\cdot)$ and $g(\cdot)$ are both
linear mappings. Recovering two variables (i.e., $\mathbf{A}$ and
$\mathbf{E}$) just from  one equation is an ill-posed problem but is still possible to be addressed  by optimizing:
\begin{equation}
\begin{array}{l}
 \text{(P0)}  ~~ \mathop {\min }\limits_{({\mathbf{A}},{\mathbf{E}})} ~~ rank(\mathbf{A}) + \lambda \left\| \mathbf{E} \right\|_{\ell _0 }  \\
 ~~~~~~~~s.t. ~~~~ \mathbf{P} = f(\mathbf{A}) +g(\mathbf{E}). \\
 \end{array}
 \label{eqs:Basic}
 \end{equation}
In (P0), $rank(\mathbf{A})$ is adopted to describe the low-rank
structure of matrix $\mathbf{A}$, and the sparse errors are
penalized via $\|\mathbf{E}\|_{\ell_0}$, where $\ell_0$ norm counts the
number of all the non-zero entries in a matrix. (P0) is always referred as sparse optimization since rank term and $\ell_0$ norm are sparse measurements for matrices and vectors, respectively. However, such sparse optimization is of little use due to the discrete nature of (P0) and the exact solution to it requires an intractable combinatorial
search.

A common approach that makes (P0) trackable tries to minimize its
convex envelope, where the rank of a matrix is replaced by the
nuclear norm and the sparse errors are penalized via $\ell_1$
norm, which are convex envelopes for $rank(\cdot)$ and $\|\cdot\|_{\ell_0}$, respectively.  Although this framework is developed on two different
norms, essentially, it is based on $\ell_1$ heuristic
since nuclear norm can be regarded as a specific case of
$\ell_1$ norm \cite{SVT-matrix-recovery}. In practical applications,  LRSL via $\ell_1$ heuristic is powerful enough for many learning tasks with relative low rank structure and sparse corruptions. However, when the desired matrix becomes complicated, e.g., it has high intrinsic rank structure or the corrupted errors become dense, the convex $\ell_1$ heuristic approaches may not  achieve promising performances. In order to handle those tough tasks via LRSL,  in this paper we take the advantages of non-convex approximation, rather than the convex $\ell_1$ heuristic, to better enhance the sparseness of signals.

We propose log-sum heuristic recovery (LHR) to  use log-sum term as the basic sparse heuristic functionality for sparse optimization in (P0). There are mainly two reasons that we use such a non-convex term to conduct LRSL. The main reason is due to its sparseness. It is indicated in \cite{reweighted} that log-sum term is a closer approximation to the $\ell_0$ norm than the $\ell_1$ norm for sparse vector representation. Therefore, it naturally inspires us to generalize its advantages from vector recovery to matrix learning.  Moreover, although the objective function derived from the log-sum term is non-convex, it is also possible to solve it by convex optimizations because the convex surrogates of log-sum term can be well defined by Taylor expansion. We will introduce an effective
non-convex optimization strategy called
\emph{majorization-minimization}(MM)
\cite{MatrixPhd}\cite{MM_surrogate} to solve it next.

MM algorithm is implemented in an iterative way that it first
replaces the non-convex component of the objective with a
convex upper-bound and then to minimize the convex upper-bound, which
exactly makes the non-convex problem fall into the general paradigm of
the reweighted schemes. Accordingly, it is possible to solve the
non-convex optimization following a sequence of convex optimizations and we will prove that with the MM framework,
 LHR  finally converges to a stationary point after successive iterations. 



LHR is a general paradigm for LRSL and we will adapt it to two specific models for practical applications. In a nutshell, LHR will be used to solve the problems of low rank matrix recovery (LRMR) and low rank representation (LRR). In LRMR,
LHR is used to recover a low rank matrix from sparse corruptions and its performance will be compared with the
benchmark Principle Component Pursuit (PCP)
\cite{Candes:1229741} method. In practice, our approach often performs
very well in spite of its simplicity. By numerical simulations,
 LHR could handle many tough tasks
that typical algorithm fails to handle. Moreover, the feasible region of LHR is much larger than PCP, which implies that it could deal with much denser corruptions and exhibits much higher rank tolerance. The feasible region of PCP subjects to the boundary of $\eta^{PCP}+\xi^{PCP}=0.35$, where $\eta$ and $\xi$ are rank rate and error rate, respectively. With the proposed LHR model, the feasible boundary can be extended to $\eta^{LHR}+\xi^{LHR}=0.58$.  The
advancements are also verified on two practical
applications of shadow removal on face images and video background
modeling.

In the second task of low rank representation, the power of LHR model will be generalized to
low rank representation for subspace clustering (SC), the goal of which aims at recovering the underlying low rank correlation of subspaces in
spite of noisy disturbances. In order to judge the performances, we will first apply it to motion segmentation in  video sequences, which is a benchmark test for SC algorithms. Besides, in order to highlight the robustness of LHR to noises and disturbances, we apply LHR to stock clustering that is to determine a stock's industrial category given its historical price record.  From both the experiments, LHR gains higher clustering accuracy than other state-of-the-art algorithms and the improvements are especially noticeable on the stock data which includes significant disturbances.

The contributions of this work are three-folds:
\begin{itemize}
  \item This work presents a log-sum heuristic recovery (LHR) algorithm to handle the typical LRSL problem with an enhanced sparsity term.  We introduce a majorization-minimization algorithm to solve the non-convex LHR optimization with reweighted schemes and theoretical justifications are provided to prove that the proposed algorithm converges to a stationary point.
  \item The proposed LHR model extends the feasible region of existing $\ell_1$ norm based LSRL algorithm, which implies that it could successfully handle  more learning tasks with  denser corruptions and higher rank.
  \item We apply the  LHR model to a new task of stock clustering which serves to demonstrate that low rank structure learning is not only a powerful tool restricted in the areas of image and vision analysis, but also can be applied to solve the profitable financial problems.
\end{itemize}

The remainder of this paper is organized as follows. We review previous works in Section \ref{sec:Relatedworks}. Section
\ref{sec:LHR} introduces the general LHR model and discusses how
to solve the non-convex LHR by MM algorithm. We addresses the low rank matrix recovery (LRMR) problem and
compare LHR model with PCP from both the simulations and
practical applications in Section
\ref{sec:REWRPCA}. The LHR model for low rank representation and subspace segmentation
is discussed in Section \ref{sec:REWLRR}. Section \ref{sec:Con} concludes this paper.

\section{Previous works}
\label{sec:Relatedworks}
In this part, we review some related works from the following
perspectives. First, we discuss two famous models in LRSL, i.e., Low Rank Matrix Recovery (LRMR) from corruptions and low rank representation (LRR).  Then, some previous works about
Majorization-Minimization algorithm and  reweighted approaches
are presented.
\subsection{Low rank structure learning}
\subsubsection{Low rank matrix recovery}
Corrupted matrix recovery considers decomposing a low rank matrix from sparse corruptions which can be formulated as $\mathbf{P}=\mathbf{A} + \mathbf{E}$,
where $\mathbf{A}$ is a low rank matrix, $\mathbf{E}$ is the sparse error and $\mathbf{P}$ is
the observed data from real world devices, e.g. cameras, sensors and
other equipments. The rank of $\mathbf{P}$ is not low, in most scenarios, due
to the disturbances of $\mathbf{E}$. How can we recover the low rank structure of
the matrix from  gross errors? This interesting topic has been discussed in a number of works,
e.g. \cite{Candes:1229741} \cite{Chandrasekaran:1183674} and
\cite{RMD}. Wright \emph{et al.} proposed the  PCP (a.k.a. RPCA)
to minimize the nuclear norm of a matrix by penalizing the
$\ell_1$ norm of errors \cite{Chandrasekaran:1183674}. PCP
could exactly recover the low rank matrix from sparse corruptions.
In some recent works, Ganesh \emph{et
al.} investigated  the parameter choosing strategy for
PCP from both the theoretical justifications and simulations
\cite{DensRPCA}. In this work, we will introduce the reweighted
schemes to further improve the performances of PCP. Our algorithm
could exactly recover a corrupted matrix from much denser errors
and higher rank.
\subsubsection{Low rank representation}
Low rank representation\cite{LRR} is a robust tool for subspace clustering \cite{SC}, the desired task of which is to classify the mixed data in their corresponding subspaces/clusters. The general model of LRR can be formulated as $\mathbf{P}=\mathbf{PA}+\mathbf{E}$, where $\mathbf{P}$ is the original mixed data, $\mathbf{A}$ is the
affine matrix that reveals the correlations between different
pairs of data and $\mathbf{E}$ is the residual of such a representation. In LRR, the affine matrix $\mathbf{A}$ is assumed to be low rank and  $\mathbf{E}$ is regarded as sparse corruptions. Compared with existing SC algorithms, LRR is much robust to noises and archives promising clustering results on public datasets.
In this work, inspired by LRR, we will introduce the log-sum recovery paradigm to LRR and show that,
with the log-sum heuristic, its robustness to corruptions can be  further improved.

 \subsection{MM algorithm and reweighted approaches }
Majorization-Minimization (MM) algorithm is widely used in  machine learning and signal processing. It is an effective strategy for non-convex problems in which the hard problem is solved  by optimizing a series of easy surrogates.  Therefore, most optimizations via MM algorithm
fall into the framework of  reweighted approaches.

In the
field of machine learning,  MM algorithm has been applied to
parameters selection for bayesian classification \cite{MMhyper}.  In the area of signal
processing, MM algorithm  leads to a number of interesting applications,
including wavelet-based processing \cite{MMwavelet} and total variation (TV)
minimization \cite{MMTV}. For compressive sensing, reweighted method
was used in $\ell_1$ heuristic and led to a number of practical
applications including  portfolio management \cite{PHDrecovery} and
image processing \cite{reweighted}. Reweighted nuclear norm was
first discussed in \cite{reweightedNuclear1} and the convergence
of such approach has been proven in \cite{Nuclear-reweighted1}.

Although there are some previous works on reweighted approaches for
rank-minimization, our approach is quite different. First, this
work tries to consider a new problem of low rank structure learning from corruptions
while not on the single task of sparse signal or nuclear
norm minimization. Besides,  existing works on reweighted nuclear norm
minimization in \cite{reweightedNuclear1}
\cite{Nuclear-reweighted1} are solved by semi-definite programming
which could only handle the matrix of relative small size. In
this paper, we will use the first order numerical algorithm (e.g., alternating direction method (ADM)) to solve
the reweighed problem, which can significantly improve the numerical performance. Due to the
distributed optimization strategy, it is possible generalize
the learning capabilities to large scale matrices.

\section{Corrupted matrix recovery via log-sum heuristic}
\label{sec:LHR}
In this part, we first discuss the widely used $\ell_1$-based method for LRSL. Then,  the log-sum heuristic (LHR) approach is proposed and we introduce the  MM algorithm to solve it. Finally, theoretical justifications are presented to prove that  LHR  can converge to a stationary point by reweighted approaches.
\subsection{$\ell_1$ heuristic for corrupted low rank matrix recovery}
As stated previously, the basic optimization (P0) is non-convex
and generally impossible to be solved as its solution usually requires
an intractable combinatorial search. In order to make
(\ref{eqs:Basic}) it trackable, convex alternatives are widely used
in a number of works, e.g. \cite{Candes:1229741}
\cite{Chandrasekaran:1183674}. Among these approaches, one
prevalent method tries to replace the rank of a matrix by its
convex envelope, i.e., the nuclear norm and the $\ell_0$ sparsity is
penalized via $\ell_1$ norm.  Accordingly, by convex relaxation, the problem
in (\ref{eqs:PCP})  can actually be recast as a semi-definite
programming.

\begin{equation}
\begin{array}{l}
~~ \mathop {\min }\limits_{({\mathbf{A}},{\mathbf{E}})} \|\mathbf{A}\|_* + \lambda \left\| \mathbf{E} \right\|_{\ell _1 }  \\
 ~~~s.t. ~~ \mathbf{P} = f(\mathbf{A}) +g(\mathbf{E}), \\
 \end{array}
 \label{eqs:PCP}
 \end{equation}
where   $\|\mathbf{A}\|_*=\sum\limits_{i = 1}^r {\sigma _i (A)}$,
is the nuclear norm of  the matrix which is defined as the
summation of the singular values of $\mathbf{A}$; and
$\|\mathbf{E}\|_{\ell_1}=\sum_{ij} |E_{ij}|$ is the $\ell_1$ norm
of a matrix. Although the objective in (\ref{eqs:PCP}) involves
two norms: nuclear norm and $\ell_1$ norm, its essence is based on
the $\ell_1$ heuristic. We will verify this point with the
following lemma.
\begin{lemma}
For a matrix $\mathbf{X}\in \mathbb{R}^{m\times n}$, its nuclear
norm is equivalent to the following optimization:
\begin{equation}
\left\| {\mathbf{X}} \right\|_*  = \left\{ \begin{array}{l}
 \mathop {\min }\limits_{({\mathbf{Y}},{\mathbf{Z}},{\mathbf{X}})} \frac{1}{2}[tr(\mathbf{Y}) + tr(\mathbf{Z})] \\
 ~~~s.t.~~\left[ {\begin{array}{*{20}c}
   {\mathbf{Y}}  \\
   {{\mathbf{X}}^T }  \\
\end{array}\begin{array}{*{20}c}
   {\mathbf{X}}  \\
   {\mathbf{Z}}  \\
\end{array}} \right] \succeq 0, \\
 \end{array} \right\}
\end{equation}
where $\mathbf{Y}\in \mathbb{R}^{m\times m},\mathbf{Z}\in
\mathbb{R}^{n\times n}$ are  both symmetric and positive definite.
The operator $tr(\cdot)$ means the trace of a matrix and $\succeq$
represents semi-positive definite. \label{Lemma:Tracenulcear}
\end{lemma}
The proof of Lemma.\ref{Lemma:Tracenulcear} may refer to
\cite{PHDrecovery}\cite{nuclear-norm}. According to this lemma, we
can replace the nuclear norm in  (\ref{eqs:PCP}) and formulate it in
the form of:
\begin{equation}
\begin{array}{l}
 \mathop {\min }\limits_{({\mathbf{Y}},{\mathbf{Z}},{\mathbf{A}},{\mathbf{E}})} \frac{1}{2}[tr({\mathbf{Y}}) + tr({\mathbf{Z}})] + \lambda \left\| {\mathbf{E}} \right\|_{\ell _1 }  \\
 ~~~s.t.~~~~~\left[ {\begin{array}{*{20}c}
   {\mathbf{Y}}  \\
   {{\mathbf{A}}^T }  \\
\end{array}\begin{array}{*{20}c}
   {\mathbf{A}}  \\
   {\mathbf{Z}}  \\
\end{array}} \right] \succeq 0 \\
 ~~~~~~~~~~~~~{\mathbf{P}} = f({\mathbf{A}}) + g({\mathbf{E}}). \\
 \end{array}
 \label{eqs:tracePCP}
\end{equation}
From Lemma.\ref{Lemma:Tracenulcear}, we know that both
$\mathbf{Y}$ and $\mathbf{Z}$ are symmetric and positive definite.
Therefore, the trace of $\mathbf{Y}$ and $\mathbf{Z}$ can be
expressed as a specific form of $\ell_1$ norm, i.e.
$tr(\mathbf{Y})=\|diag(\mathbf{Y})\|_{\ell_1}$. $diag(\mathbf{Y})$ is an operator that
only keeps the entries on the diagonal position of $\mathbf{Y}$ in
a vector. Therefore, the optimization in (\ref{eqs:tracePCP}) can
be expressed as:
\begin{equation}
 \mathop {\min }\limits_{\hat{X}\in \hat{D}} \frac{1}{2}(\|diag({\mathbf{Y}})\|_{\ell_1} + \|diag({\mathbf{Z}})\|_{\ell_1}) + \lambda \left\| {\mathbf{E}} \right\|_{\ell _1 }, \\
 \label{eqs:l1H}
\end{equation}
where $\hat{X}=\{\mathbf{Y,Z,A,E}\}$ and
\[
\hat{D} = \{ ({\bf{Y}},{\bf{Z}},{\bf{A}},{\bf{E}}):\left[ {\begin{array}{*{20}c}
   {\mathbf{Y}}  \\
   {{\mathbf{A}}^T }  \\
\end{array}\begin{array}{*{20}c}
    {\mathbf{A}}  \\
   {\mathbf{Z}}  \\
\end{array}} \right] \succeq 0,({\bf{A}},{\bf{E}}) \in C\}.
\]
$({\bf{A}},{\bf{E}}) \in C$ stands for convex constraint.
\subsection{Log-sum heuristic for matrix recovery}
By Lemma.\ref{Lemma:Tracenulcear}, the convex problem with
two norms in (\ref{eqs:PCP}) has been successfully converted to
an optimization only with $\ell_1$ norm and therefore it is  called
$\ell_1$-heuristic.
 $\ell_1$ norm is the convex envelope of the concave $\ell_0$ norm but a
number of previous research works have indicated the limitation of
approximating $\ell_0$ sparsity with $\ell_1$ norm, e.g.,
\cite{reweighted}\cite{elastic-net}. It is natural to ask, for
example, whether might a different alternative not only find a correct solution, but also outperform the performance of
$\ell_1$ norm?  Next we will introduce the log-sum term
to represent the sparsity of signals.
\begin{definition}
For any matrix $\mathbf{X}\in \mathbb{R}^{m\times n}$, the log-sum term is defined as $\|\mathbf{X}\|_L=\sum_{ij} \log(|X_{ij}|+\delta)$, where $\delta>0$ is a small regularization constant.
\end{definition}

The prominent reason that we use this term is mainly due to its sparsity. As indicated in \cite{reweighted},  the log-sum term lies between the scope of the $\ell_0$ norm and
$\ell_1$ norm, which makes it be a closer
approximation of $\ell_0$  norm
\cite{reweighted}\cite{PHDrecovery} and therefore, it is used to encourage the sparsity in the optimization.
We propose  Log-sum Heuristic Recovery (LHR) model $H(\hat{X})$:
\begin{equation}
\begin{array}{l}
\text{(LHR)} H(\hat{X})= \mathop {\min }\limits_{\hat{X}\in {\hat {D}}} \frac{1}{2}(\|diag({\mathbf{Y}})\|_{L} + \|diag({\mathbf{Z}})\|_{L}) + \lambda \left\| {\mathbf{E}} \right\|_{L}.  \\
 \end{array}
 \label{eqs:LHR}
\end{equation}

From the formulation of LHR, obviously, it differs from (\ref{eqs:l1H}) only on the selection of
the sparse norm, where the later uses log-sum term instead of the
typical $\ell_1$ norm. Although we have placed a powerful term to
enhance the sparsity in LHR model, unfortunately, it also causes
non-convexity into the objective function. The LHR model is not
convex since the log-function over $\mathbb{R_{++}}=(\delta,\infty)$ is
concave. In most cases, non-convex problem can be extremely
hard to solve. Fortunately, the convex upper bound of
$\|\cdot\|_L $ can be easily found and defined.
Therefore, we will introduce the majorization-minimization
algorithm to solve the LHR optimization.

\subsection{The majorization-minimization for LHR optimization}
The majorization-minimization (MM) algorithm replaces the hard problem
by a sequence of easier ones. It proceeds in an Expectation
Maximization (EM)-like fashion by  repeating  two steps of
\textbf{M}ajorization and \textbf{M}inimization in an iterative
way. During the \emph{Majorization} step, it constructs the convex
upper bound of the non-convex objective. In the
\emph{Minimization} step, it minimizes the upper bound.

To see how the MM works for LHR, Let's recall the objective function in (\ref{eqs:LHR}) and make some simple algebra operations:
\begin{equation}
\begin{array}{l}
 \frac{1}{2}[\left\| {diag({\bf{Y}})} \right\|_L  + \left\| {diag({\bf{Z}})} \right\|_L ] + \lambda \left\| {\bf{E}} \right\|_L  \\
  = \frac{1}{2}[\sum\nolimits_i {\log (Y_{ii}+\delta )}  + \sum\nolimits_k {\log (Z_{kk}+\delta )} ]
  \\~~~+ \lambda \sum\nolimits_{ij} {\log (|E_{ij}|+\delta )}  \\
  = \frac{1}{2}[\log \det ({\bf{Y}+\delta \mathbf{I}_{m}}) + \log \det({\bf{Z}+\delta \mathbf{I}_{n}})] 
 \\~~~+ \lambda \sum\nolimits_{ij} {\log (|E_{ij}|+\delta )}  \\
 \end{array}
 \label{eqs:expasion}
\end{equation}
where $\mathbf{I}_m \in \mathbb{R}^{m \times m}$ is an identity matrix. It is well known that the concave function is bounded by its
first-order  Taylor expansion. Therefore, we calculate the convex
upper bounds of all the terms in (\ref{eqs:expasion}). For the
term $\log \det(\mathbf{Y}+\delta\mathbf{I}_m)$,
\begin{equation}
\begin{array}{l}
 \log \det ({\bf{Y}+\delta \mathbf{I}_{m}}) \le \log \det ({\bf{\Gamma}}_{Y}+ \delta \mathbf{I}_{m})  \\~~~~~~~~~~~~~~~~~~~~~~~~+tr[({\bf{\Gamma}}_Y +\delta \mathbf{I}_{m} )^{ - 1}(\bf{Y} - {\bf{\Gamma}}_Y )]. \\
 \end{array}
 \label{eqs:MatrixExpasion}
\end{equation}
The inequality in (\ref{eqs:MatrixExpasion}) holds for any $\mathbf{\Gamma}_Y\succ 0$. Similarly, for any $(\Gamma_E)_{ij}>0$,
\begin{equation}
\begin{array}{l}
\sum\nolimits_{ij} {\log (|E_{ij}|+\delta)}  \le \sum\nolimits_{ij} {[\log [(\Gamma_E)_{ij}+\delta]+ \frac{E_{ij} - (\Gamma_E)_{ij} }{{(\Gamma_E)_{ij}+\delta }}]}
 \end{array}
 \label{eqs:ErroExpansion}
\end{equation}

We replace each term in (\ref{eqs:expasion}) with the convex upper
bound and define $T(\hat{X}|\hat{\Gamma})$ as the
surrogate function after convex relaxation. Therefore, we can
instead optimize the following problem
\begin{equation}
\begin{array}{l}
\mathop{\min}\limits_{\hat{X}\in {\hat {D}}} T(\hat{X}|\hat{\Gamma})=%
\frac{1}{2}tr[({\bf{\Gamma }}_Y  + \delta {\bf{I}}_m )^{ - 1} {\bf{Y}}] + \frac{1}{2}tr[({\bf{\Gamma }}_z  + \delta  {\bf{I}}_n )^{ - 1} {\bf{Z}}]\\~~~~~~~~~~~~~~~~~~ + \lambda \sum\nolimits_{ij} {(\Gamma _{E_{ij} }  + \delta  )^{ - 1} E_{ij} }  + const,
\end{array}
\label{eqs:T}
\end{equation}
In (\ref{eqs:T}),  set $\hat{X}=\{\mathbf{Y},\mathbf{Z},\mathbf{A},\mathbf{E}\}$, which contains all the variables to
be optimized and set $\hat{\Gamma}=\{\mathbf{\Gamma}_Y, \mathbf{\Gamma}_Z,
\mathbf{\Gamma}_E\}$ contains all the parameter matrices. At the end of (\ref{eqs:T}), \emph{const} stands for the constants that are irrelative to $\{\mathbf{Y},\mathbf{Z},\mathbf{A},\mathbf{E}\}$. In  some
previous works of MM algorithms \cite{MMhyper} \cite{reweighted}
\cite{Nuclear-reweighted1}, they denote the parameter  $\hat{\Gamma}$ in
$t^{th}$ iteration with the optimal value of $\hat{X}$ of the
last iteration, i.e. $\hat{\Gamma}=\hat{X^t}^*$. According to the
discussions above, we provide the MM algorithm for (LHR)
minimization in Algorithm \ref{alg:MM-LHR}.
\begin{algorithm}
\SetKwInOut{Majorization}{Majorization}\SetKwInOut{Minimization}{Minimization}
\SetKwData{Define}{Define}
\SetKwInOut{Initialization}{Initialization}

\Initialization{$t:= 0$ }
 \Repeat { convergence}
 { \Majorization\
   $\hat{\Gamma}^{t}:=\hat{X}^{t}$\;
   \Define convex upper bound $T(\hat{X}|\hat{\Gamma}^{t})$\;
 \Minimization\
 $\hat{X}^{t+1}=\arg \mathop{\min.}\limits_{\hat{X}\in {\hat {D}}} T(\hat{X}|\hat{\Gamma}^{t})$\;\label{line:mg}
$t:=t+1$\;}
\caption{A MM algorithm for LHR minimization}\label{alg:MM-LHR}
\end{algorithm}
Before elaborately discussing how to numerically solve the
optimization, we will first discuss some theoretical properties of
it.
\subsection{Theoretical justifications}
In this part, for simplicity,
we define the objective function in (\ref{eqs:LHR}) as $H(\hat{X})$
and the surrogate function in (\ref{eqs:T})
is defined as $T(\hat{X}|\hat{\Gamma})$. $\hat{X}$ is a set containing
all the variables  and  set $\hat{\Gamma}$ records the parameter
matrices.  The convergence property of general MM algorithm was separately distributed on some early mathematical journals \cite{hunter2005variable}\cite{lange1995gradient} which is a bit obscure and were not generally read by researchers in the community of computer science. Besides, previous works on MM convergence are almost on the variable selection models. In this paper,  we  specify it to our LHR model and try to explain it in a plain way. Before discussing the
convergence property of LHR, we will first provide two lemmas.
\begin{lemma}
If set $\hat{\Gamma}^{t}:=\hat{X}^{t}$, MM algorithm could monotonically decrease the non-convex objective function $H(\hat{X})$, i.e. $H(\hat{X}^{t+1}) \leq H(\hat{X}^{t})$.
    \label{Lemma:Decrease}
\end{lemma}

\begin{proof}
 In order to prove the monotonically decrease property, we can instead prove:
\begin{equation}\label{eqs:MDproof}
    \begin{array}{l}
    H(\hat{X}^{t+1}) \leq T(\hat{X}^{t+1}|\hat{\Gamma}^{t})\leq T(\hat{X}^{t}|\hat{\Gamma}^{t})
    =H(\hat{X}^{t}).
    \end{array}
\end{equation}
We prove (\ref{eqs:MDproof}) by the following three steps:

(\textbf{i}) The first inequality follows from the
argument  that $T(\hat{X}|\hat{\Gamma})$ is the upper-bound of
$H(\hat{X})$.

(\textbf{ii}) The second inequality holds since the MM algorithm
computes  $\hat{X}^{t+1}=\arg \mathop {\min }\limits_{\hat{X}} T(\hat{X}|\hat{\Gamma}^t)$. The function $T(\cdot)$ is convex,
therefore, $\hat{X}^{t+1}$ is the unique global minimum. This property
guarantees that $T(\hat{X}^{t+1}|\hat{\Gamma}^{t+1})<
T(\cdot|\hat{\Gamma}^{t})$  with any
$\hat{X} \neq \hat{X}^{t+1}$ and $T(\hat{X}^{t+1}|\hat{\Gamma}^{t+1})=
T(\cdot|\hat{\Gamma}^{t})$  if and only if
$\hat{X} =\hat{X}^{t+1}$.

(\textbf{iii}) The last equality can be easily verified by expanding
$T(\hat{X}^{t}|\hat{\Gamma}^{t})$ and making some simple algebra.
The transformation is straightforward and omitted here.
\end{proof}

\begin{lemma}
Let $\hat{X}=\{\hat{X}^0,\hat{X}^1...\hat{X}^t...\}$ be a sequence generated by MM framework in Algorithm \ref{alg:MM-LHR}, after successive iterations, such a sequence converges to the same limit point.
\label{Lemma:Convergence}
\end{lemma}
\begin{proof}
  We give a proof by contradiction. We assume that sequence $\hat{X}$ diverges, which means that $\lim\limits_{t\rightarrow\infty}\|\hat{X}^{t+1}-\hat{X}^{t}\|_F\ne0$. According to the discussions in Appendix \ref{sec:subconvergent}, we know that there exists a convergent subsequence $\hat{X}^{t_k}$  converging to $\phi$, i.e. $\lim\limits_{k\rightarrow\infty}\hat{X}^{t_k}=\phi$ and meanwhile, we can construct another convergent subsequence $\hat{X}^{t_k+1}$ that $\lim\limits_{k\rightarrow\infty}\hat{X}^{t_k+1}=\varphi$. We assume that $\phi\neq \varphi$. Since the convex upper-bound $T(\cdot|\hat{\Gamma})$ is continuous, we get
$\mathop {\lim }\limits_{k \to \infty } T(\hat X^{t_k + 1} |\hat{\Gamma}^{t_k} ) = T(\underbrace {\mathop {\lim }\limits_{k \to \infty } \hat X^{t_k + 1} }_\varphi |\hat{\Gamma}^{t_k} ) < T(\underbrace {\mathop {\lim }\limits_{k \to \infty } \hat X^{t_k} }_\phi |\hat{\Gamma}^{t_k} ) = \mathop {\lim }\limits_{k \to \infty } T(\hat X^{t_k} |\hat{\Gamma}^{t_k} )$
The \emph{strict less than operator}  "$<$" holds because $\varphi \ne \phi$. See  $(\textbf{ii})$ in the proof of Lemma \ref{Lemma:Decrease} for details.
Therefore, it is straightforward to get the following inequalities:
$\mathop {\lim }\limits_{k \rightarrow \infty } H(\hat X^{t_k + 1} ) \le \mathop {\lim }\limits_{k \rightarrow \infty } T(\hat X^{t_k + 1} |\hat{\Gamma}^{t_k} ) < \mathop {\lim }\limits_{k \rightarrow \infty } T(\hat X^{t_k} |\hat{\Gamma}^{t_k} ) = \mathop {\lim }\limits_{k \rightarrow \infty } H(\hat X^{t_k} )$. Accordingly,
\begin{equation}
\mathop {\lim }\limits_{k \to \infty } H(\hat X^{t_k + 1} )< \mathop {\lim }\limits_{k \to \infty } H(\hat X^{t_k} )
\label{eqs:1st}
\end{equation}
Besides, it is obvious that the function of $H(\cdot)$ in (\ref{eqs:LHR}) is bounded below, i.e. $H(\hat{X})>(mn+m+n)\log\delta$. Moreover, as proved in Lemma \ref{Lemma:Decrease},  $H(\hat{X})$ is monotonically decreasing, which guarantee that
$\mathop {\lim }\limits_{t \to \infty } H(\hat X^{t } )$
exists, i.e.
\begin{equation}
\begin{array}{l}
\mathop {\lim }\limits_{k \to \infty } H(\hat X^{t_k } )=\mathop {\lim }\limits_{t \to \infty } H(\hat X^{t} )=\mathop {\lim }\limits_{t \to \infty } H(\hat X^{t+1} )\\~=\mathop {\lim }\limits_{k \to \infty } H(\hat X^{t_k + 1} )
\end{array}
\label{eqs:2nd}
\end{equation}
Obviously, (\ref{eqs:2nd}) contradicts to (\ref{eqs:1st}). Therefore, $\phi=\varphi$ and we get the conclusion that $\lim\limits_{t\rightarrow\infty}\|\hat{X}^{t+1}-\hat{X}^{t}\|_F=0$.
\end{proof}

Based on the two lemmas proved previously, we can give the convergence theorem of the proposed LHR model.
\begin{theorem}
With the MM framework, LHR model  finally converges to a stationary point.
\end{theorem}

\begin{proof}
  As stated in Lemma \ref{Lemma:Convergence}, the sequences generated by MM algorithm converges to a limitation and here we will first prove that the convergence is a \emph{fixed point}. We define the mapping from $\hat{X}^k$ to $\hat{X}^{k+1}$ as $M(\cdot)$, and it is straightforward to get,
$\mathop {\lim }\limits_{t \to \infty } \hat X^t  = \mathop {\lim }\limits_{t \to \infty } \hat X^{t + 1}  = \mathop {\lim }\limits_{t \to \infty } M(\hat X^t )$, which implies that
$\mathop {\lim }\limits_{t \to \infty } \hat X^t=\phi$ is a fixed point. In the MM algorithm, when constructing the upper-bound, we use the first-order Taylor expansion.  It is well known that the convex surrogate $T(\hat{X}|\hat{\Gamma})$ is tangent to $H(\hat{X})$ at $\hat{X}$ by the property of Taylor expansion. Accordingly, the gradient vector of $T(\hat{X}|\hat{\Gamma})$ and $H(\hat{X})$ are equal when evaluating at $\hat{X}$. Besides, we know that at the fixed point,
${\bf{0}} \in \nabla _{\hat X = \phi } T(\hat X{\rm{|\hat{\Gamma}}})$ and because it is tangent to $H(\hat{X})$, we can directly get that ${\bf{0}} \in \nabla _{\hat X = \phi } H(\hat{X})$ which proves that the convergent fixed point $\phi$ is also a \emph{stationary point} of $H(\cdot)$.
\end{proof}

In this part, we have shown that with the MM algorithm, LHR model could converge to a stationary point. However, it is impossible to claim that the converged point is the global minimum since the objective function of LHR is not convex. Fortunately, with a good starting point, we can always find a desirable solution by iterative approaches. In this paper, the solution of $\ell_1$ heuristic model was used as a starting point and it could always lead to a satisfactory result.
\subsection{Solve LHR via reweighted approaches}
To numerically solve the LHR optimization, we remove the
constants that are irrelative to $\mathbf{Y,Z}$ and $\mathbf{E}$ in $T(\hat{X}|\hat{\Gamma})$ and get the new convex objective
\[
 \min \frac{1}{2}[tr({\bf{W}}_Y^2  {\bf{Y}}) + tr({\bf{W}}_Z^2  {\bf{Z}})] + \lambda \sum\nolimits_{ij} {(W_E )_{ij} } E_{ij}  \\
\]
where $\bf{W}_{Y(Z)}=(\bf{\Gamma}_{Y(Z)}+\delta \bf{I}_{m(n)})^{-1/2}$ and $(W_E)_{ij}=(E_{ij}+\delta)^{-1},\forall ij$.  It is worth noting that $tr({\bf{W}}_Y^2  {\bf{Y}})=tr({\bf{W}}_Y  {\bf{Y}} \bf{W}_Y )$. Besides, since both $\mathbf{W}_Y$ and $\mathbf{W}_Z$ are positive definite, the first constraint in (\ref{eqs:LHR})
 is equivalent to
\[
\left[ {\begin{array}{*{20}c}
   {{\bf{W}}_Y } & {\bf{0}}  \\
   {\bf{0}} & {{\bf{W}}_Z }  \\
\end{array}} \right]\left[ {\begin{array}{*{20}c}
   {\bf{Y}} & {\bf{A}}  \\
   {{\bf{A}}^T } & {\bf{Z}}  \\
\end{array}} \right]\left[ {\begin{array}{*{20}c}
   {{\bf{W}}_Y } & {\bf{0}}  \\
   {\bf{0}} & {{\bf{W}}_Z }  \\
\end{array}} \right] \succeq {\bf{0}}
\]

Therefore, after convex relaxation, the optimization in (\ref{eqs:LHR}) now subjects to
\begin{equation}
\begin{array}{l}
 \min ~~\frac{1}{2}[tr({\bf{W}}_Y  {\bf{YW}}_Y ) + tr({\bf{W}}_Z  {\bf{ZW}}_Z  )] + \lambda \|{\bf{W}}_{E}  \odot {\bf{E}}\|_{\ell_1}  \\
 s.t.~~~
\left[ {\begin{array}{*{20}c}
   {{\bf{W}}_Y {\bf{YW}}_Y } & {{\bf{W}}_Y {\bf{{\rm A}W}}_Z }  \\
   {({\bf{W}}_Y {\bf{{\rm A}W}}_Z )^T } & {{\bf{W}}_Z {\bf{ZW}}_Z }  \\
\end{array}} \right] \succeq {\bf{0}}

\\
 ~~~~~~~~~{\bf{P}} = f({\bf{A}}) + g({\bf{E}}) \\
 \end{array}
 \label{eqs:majorization}
\end{equation}

Here, we apply Lemma \ref{Lemma:Tracenulcear} to (\ref{eqs:majorization}) once again and rewrite the optimization in (\ref{eqs:majorization}) in the form of the summation of the nuclear norm and $\ell_1$ norm,
\begin{equation}
\begin{array}{l}
 \mathop {\min }\limits_{({\bf{A}},{\bf{E}})}. \left\| {\bf{W}}_Y {\bf{AW}}_Z \right\|_*  + \lambda \|{\bf{W}}_{E}  \odot {\bf{E}}\|_{\ell_1} \\
~~ s.t.~{\bf{P}} = f({\bf{A}}) + g({\bf{E}}) \\
 \end{array}
\label{eqs:rew}
\end{equation}

In (\ref{eqs:rew}), the operator $\odot$ in the error term denotes the
component-wise product of two variables, i.e., for $\mathbf{W}_{E}$ and $\mathbf{E}$:
$(\mathbf{W}_{E}\odot \mathbf{E})_{ij}=(W_{E})_{ij}E_{ij}$. According to \cite{PHDrecovery}, we know that $\mathbf{Y}^*=\mathbf{U}\Sigma\mathbf{U}^T$ and $\mathbf{Z}^*=\mathbf{V}\Sigma\mathbf{V}^T$, if we do singular value decomposition for $\mathbf{A}^*=\mathbf{U}\Sigma\mathbf{V}^T$. Accordingly, the weight matrix $\mathbf{W}_Y=(\mathbf{U}\Sigma\mathbf{U}^T+\delta \mathbf{I}_m)^{-1/2}$ and matrix $\mathbf{W}_{Z}=(\mathbf{V}\Sigma\mathbf{V}^T+\delta \mathbf{I}_n)^{-1/2}$.

 Here, based on MM algorithm, we have converted the non-convex LHR optimization to be  a sequence of convex reweighted problems. We call it \emph{reweighted method} (\ref{eqs:rew}) since in each iteration we should re-denote the weight matrix set $\hat{W}$ and use the updated weights to construct the surrogate convex function. Besides, the objective  in (\ref{eqs:rew}) is convex with a summation of a nuclear norm and a $\ell_1$ norm and can be solved by convex optimization. In the next two sections, the general LHR model will be adapted to two specific models and we will provide the optimization strategies for those two models,respectively.

\section{Low rank matrix recovery from corruptions}
\label{sec:REWRPCA}
In this part, we first apply the LHR model to recover a low rank matrix from corruption and its performance is compared with the widely used Principal Component Pursuit (PCP).
\subsection{Joint optimization for LHR}
\label{sec:RPCA-solve}
Based on the LHR derivations, the corrupted low rank matrix recovery problem can be formulated as a reweighted problem:
\begin{equation}
\begin{array}{l}
 \mathop {\min }\limits_{({\bf{A}},{\bf{E}})}. \left\| {{\bf{W}}_Y  {\bf{AW}}_Z  } \right\|_*  + \lambda \|{\bf{W}}_{E}  \odot {\bf{E}}\|_{\ell_1} \\
~~ s.t.~~{\bf{P}} = \bf{A} + \bf{E} \\
 \end{array}
\label{eqs:rewRPCA}
\end{equation}
Due to the reweighted weights are placed in the nuclear norm, it is impossible to directly get the closed-form solution of the nuclear norm minimization. Therefore, inspired by the work \cite{LRR}, we  introduce another variable $\mathbf{J}$  to (\ref{eqs:rewRPCA}) by adding another equality constraint and to solve,
\begin{equation}
\begin{array}{l}
 \min .\left\| {\bf{J}} \right\|_*  + \lambda \left\| {{\bf{W}}_{E }  \odot {\bf{E}}} \right\|_{\ell _1 }  \\
 s.t.~~{\bf{h}}_1  = {\bf{P}} - {\bf{A}} - {\bf{E}} = \bf{0} \\
 ~~~~~~{\bf{h}}_2  = {\bf{J}} - {\bf{W}}_Y {\bf{AW}}_Z  = \bf{0}\\
 \end{array}
 \label{eqs:REWRPCAE}
\end{equation}
Based on the transformation, there is only one single $\bf{J}$ in the nuclear norm of the objective that we can directly get its closed-form update rule by \cite{SVT-matrix-recovery}. There
are quite a number of methods that can be used to solve it, e.g.
with  Proximal Gradient (PG)
 algorithm \cite{APA}
or Alternating Direction Methods (ADM) \cite{ADM}. In this paper, we
will introduce the ADM method since it is more effective and
efficient. Using the ALM method \cite{ALM}, it is computationally expedient to relax the equality in (\ref{eqs:REWRPCAE}) and instead solve:
\begin{equation}
\begin{array}{l}
 L = \left\| {\bf{J}} \right\|_*  + \lambda \left\| {{\bf{W}}_{E }  \odot {\bf{E}}} \right\|_{\ell _1 }  +  < {\bf{C}}_1 ,{\bf{h}}_1  >  \\
  +  < {\bf{C}}_2 ,{\bf{h}}_2  >  + \frac{\mu }{2}(\left\| {{\bf{h}}_1 } \right\|_F^2  + \left\| {{\bf{h}}_2 } \right\|_F^2 ) \\
 \end{array}
 \label{eqs:ALM-RRPCA}
\end{equation}
where $<,>$ is an inner product and $\mathbf{C}_1$ and $\mathbf{C}_2$ are the lagrange multipliers, which can be updated via dual ascending method.
(\ref{eqs:ALM-RRPCA}) contains three variables, i.e., $\mathbf{J},\mathbf{E}$ and  $\mathbf{A}$. Accordingly, it is possible to solve problem via a distributed optimization strategy called Alternating Direction Method (ADM).  The convergence of the ADM for convex problems has been widely discussed  in a number of works \cite{ADM}\cite{ADM_covergence}. By ADM, the joint optimization can be  minimized  by four steps as $\mathbf{E}$-minimization, $\mathbf{J}$-minimization, $\mathbf{A}$-minimization and dual ascending. We first provide the update rule for $\mathbf{E}$-minimization,
\begin{equation}
{\bf{E}} = \mathop {\arg \min }\limits_{\bf{E}} \lambda \left\| {{\bf{W}}_{E}  \odot {\bf{E}}} \right\|_*  + \frac{\mu }{2}\left\| {{\bf{E}} - ({\bf{P}} - {\bf{A}} + \mu ^{ - 1} {\bf{C}}_2 )} \right\|_{_F }^2
\end{equation}
It is well known (see, for example,
\cite{Donoho-denoising-by-threshold}) that for scalars $x$ and $y$,
the unique optimal solution to the problem
\begin{equation}
\mathop {\min }\limits_x \alpha \left| x \right| + \frac{1}{2}(x -
y)^2 \label{eqs:l1}
\end{equation}
is given by
\begin{equation}
x^* {\rm{ }} = {\rm{ }}{\mathop{\rm sgn}} (y){\rm{ }}\max (|y| -
\alpha ,{\rm{ }}0){\rm{ }} \buildrel\textstyle.\over= {\rm{
}}s_\alpha  (y).
\end{equation}
$\mathbf{E}^*$ is a solution to the $\mathbf{E}$-minimization if and only if for all $i,
j$,
\begin{equation}
E_{ij}^*  = s_{\lambda \mu ^{ - 1} \left| {W_{ij} } \right|} (P - A
- \mu ^{ - 1} C)_{ij}
\end{equation}
Similarly, $\mathbf{J}$-minimization can be solved by
\begin{equation}
{\bf{J}} = \mathop {\arg \min }\limits_{\bf{J}} \left\| {\bf{J}} \right\|_*  + \frac{\mu }{2}\left\| {{\bf{J}} - {\bf{W}}_{Y} {\bf{AW}}_{Z}  + \mu ^{ - 1} {\bf{C}}_1 } \right\|_{_F }^2
\label{eqs:J-min}
\end{equation}
For matrices $\mathbf{X},\mathbf{D}$, previous works, e.g. \cite{nuclear-norm}
\cite{CandesMatrxiRecoveryViaConvex}, have indicted that the unique closed-form
optimal solution to the problem
\begin{equation}
\mathop {\min }\limits_{\mathbf{X}} \alpha \left\| \mathbf{X} \right\|_*  +
\frac{1}{2}\left\| {\mathbf{X} - \mathbf{D}} \right\|_F^2
\end{equation}
is given by
\begin{equation}
\mathbf{X}^*  = \mathbf{U}s_\alpha  (\Sigma )\mathbf{V}^T  \buildrel\textstyle.\over= d_\alpha
(\mathbf{D}), \label{eqs:MS}
\end{equation}
where $\mathbf{D}=\mathbf{U}\Sigma \mathbf{V}^T$ denotes the singular value decomposition of
$\mathbf{D}$. From (\ref{eqs:MS}), it is immediate that the unique optimal
solution to (\ref{eqs:J-min}) is given by
\begin{equation}
\mathbf{J}^*  = d_{\mu^{-1}}  ({\bf{W}}_{Y} {\bf{AW}}_{Z}  + \mu ^{ - 1} {\bf{C}}_2 ).
\end{equation}
Finally, the solution to $\mathbf{A}$ is based on the following optimization problem,
\begin{equation}
{\bf{A}^*} = \mathop {\arg \min }\limits_{\bf{A}} \left\| {{\bf{h}}_1  + \mu ^{ - 1} {\bf{C}}_1 } \right\|_{_F }^2  + \left\| {{\bf{h}}_2  + \mu ^{ - 1} {\bf{C}}_2 } \right\|_{_F }^2
\end{equation}
which is only a summation of two $F$-norms that can be addressed by gradient-descending method. Here, we provide the update-rule that
\[
{\bf{A}}^{k + 1}  = {\bf{A}}^k  +\gamma [  {\bf{W}}_Y ({\bf{h}^k}_1  + \mu ^{ - 1} {\bf{C}}_1 ){\bf{W}}_Z  + ({\bf{h}}_2^k  + \mu ^{ - 1} {\bf{C}}_2 )]
\]
Based on all the previous discussions, we can now give the whole framework to solve the LHR model for LRMR via reweighted schemes in Algorithm \ref{alg:REW}.

\begin{algorithm}
\SetKwInOut{Majorization}{Majorization}\SetKwInOut{Minimization}{Minimization}
\SetKwData{Reset}{Reset}
\SetKwInOut{Initialization}{Initialization}\SetKwInOut{Input}{Input}\SetKwInOut{Output}{Output}

\Input{Corrupted matrix $P$ and parameter $\lambda$}
\Initialization{$t:= 1 ,E^{0}_{ij}:=1, \forall i,j.~ \mathbf{W}_{Y(Z)}^{(1)}=\mathbf{I}_{m(n)}$. }
 \Repeat {convergence}
    {\tcp{Dynamically update the weight matrix.}
    $ \mathbf{W}_E^{(t)}:=(|\mathbf{E}^{(t-1)}|+ \delta_1)^{-1}$ \;
    $\mathbf{U}\Sigma\mathbf{V}^T=SVD(\mathbf{A}^{(t-1)})$\;
    $ \mathbf{W}_Y^{(t)}:=(\mathbf{U}\Sigma\mathbf{U}^T+ \delta_2\mathbf{I}_m)^{-1/2}$ \;
    $ \mathbf{W}_Z^{(t)}:=(\mathbf{V}\Sigma\mathbf{V}^T+ \delta_2\mathbf{I}_n)^{-1/2}$ \;
      \Reset $C_{0}>0;\mu_0>0;\rho>1;k=1;\mathbf{A}^0=\mathbf{E}^0=\bf{0}$\; \label{line:REWstr}
      \While{not converged}
          {      \tcp{Variables updating.}
          $E_{ij}^{k}  = s_{\lambda \mu ^{ - 1} \left| {(W_E^{(t)})_{ij} } \right|} (P - A^{k-1}
- \mu ^{ - 1} C_1^k)_{ij},\forall ij$\;
          $\mathbf{J}^{k}  = d_{\mu^{-1}}  ({\bf{W}}_{Y}^{(t)} {\bf{A}^{k-1}\mathbf{W}}^{(t)}_{Z}  + \mu ^{ - 1} {\bf{C}}_2^k )$\;
          ${\bf{A}^k}  = {\bf{A}}^{k-1}  + \gamma [ -{\bf{W}}_Y^{(t)} ({\bf{h}}_1^k  + \mu ^{ - 1} {\bf{C}}_2^k ){\bf{W}}_Z^{(t)}
          +({\bf{h}}_2^k  + \mu ^{ - 1} {\bf{C}}_1^k )]$\;
          \tcp{Dual ascending.}
          $\mathbf{C}_1^k=\mathbf{C}_1^{k-1}+\mu_k \mathbf{h}_1^k$\;
          $\mathbf{C}_2^k=\mathbf{C}_2^{k-1}+\mu_k \mathbf{h}_2^k$\;
          $k:=k+1, \mu_{k+1}=\rho \mu_k$;
          } \label{line:REWend}
$(\mathbf{A}^{(t)},\mathbf{E}^{(t)})=(\mathbf{A}_k,\mathbf{E}_k)$\;
$t:=t+1$\;
}
\Output{$(\mathbf{A}^{(t)},\mathbf{E}^{(t)})$.}
\caption{Optimization strategy of LHR for corrupted matrix recovery}
\label{alg:REW}
\end{algorithm}

\subsection{Numerical simulations}
We have explained how to recover a low rank matrix via LHR in preceding sections. In this section, we will conduct some experiments to test
its performances
with the comparisons to  robust PCP from both the simulations and practical data.
\subsubsection{General evaluation}
We demonstrate the accuracy of the proposed LHR algorithm on
randomly generated matrices. For an equivalent
comparison, we adopted the same data generating method
in \cite{Candes:1229741} that all the algorithms are
performed on the squared matrices and the ground-truth low rank matrix (rank $r$)
with $m \times n$ entries, denoted as $\mathbf{A}^*$ , is generated by
independent random orthogonal model \cite{Candes:1229741}; the sparse
error $\mathbf{E}^*$ is generated via uniformly sampling the matrix and the
error values are randomly generate in the range [-100,100].
The corrupted matrix is generated by
$\mathbf{P}=\mathbf{A}^*+\mathbf{E}^*$, where $\mathbf{A}^*$ and $\mathbf{E}^*$ are the ground truth. For simplicity,
we denote the rank rate as $\eta=\frac{rank(\mathbf{A}^*)}{max\{m,n\}}$ and the error rate
as $\xi=\frac{\|\mathbf{E}\|_{\ell_0}}{m\times n}$.

For an equivalent comparison, we use the code in \cite{exact_alm} to solve the PCP problem \footnote{ In \cite{ALM}, Lin \emph{et al.} provided two solvers, i.e. exact and inexact solvers, to solve the PCP problem. In this paper, we use the exact solver for PCP because it performs better than inexact solver.}. \cite{Candes:1229741} indicated that PCP method could exactly recover a low rank matrix from corruptions within the region of $\eta+\xi < 0.35$. Here, in order to highlight the effectiveness of our LHR model, we directly consider much difficult tasks that we set $\eta+\xi=0.5$. Each experiment is repeated for ten times and the median values
\footnote{We do not use the average values here since in cases of
divergence some extreme large outliers may greatly affect the
average values of the accuracy.} are tabulated in
Table.\ref{tab:evaluation}. In the table,
$\frac{\|\mathbf{A}-\mathbf{A}^*\|_F}{\|\mathbf{A}^*\|_F}$ denotes the recovery accuracy, $rank$ denotes the rank of the
recovered matrix $\mathbf{A}$, $\|\mathbf{E}\|_{\ell_0}$ is the card of the recovered
errors and $time$ records the computational costs (in seconds).
\begin{table*}
\caption{Evaluations of low-rank matrix recovery  of Robust PCA and
Log-sum Heuristic Recovery. } \centering
\begin{tabular}{|c|c|cccc|cccc|}
\hline
           &            & \multicolumn{ 4}{|c}{$rank(\mathbf{A}^*)=0.4m$~~~$||\mathbf{E}^*||_{\ell_0}=0.1m^2$    } & \multicolumn{ 4}{|c|}{$rank(\mathbf{A}^*)=0.1m$~~~$||\mathbf{E}^*||_{\ell_0}=0.4m^2$    } \\
\hline
      $m=n$ &    methods &    $\frac{\|\mathbf{A}-\mathbf{A}^*\|_F}{\|\mathbf{A}^*\|_F}$  &    $rank(\mathbf{A})$ & $||\mathbf{E}||_{\ell_0}$  &    $time(s)$ &    $\frac{\|\mathbf{A}-\mathbf{A}^*\|_F}{\|\mathbf{A}^*\|_F}$  &    $rank(\mathbf{A})$ & $||\mathbf{E}||_{\ell_0}$  &    $time(s)$ \\
\hline \hline
\multicolumn{ 1}{|c|}{200} &   PCP &   4.6e-1 &         103 &       21066 &        5.9 &   1.2e-1 &        107 &      23098 &        7.4 \\

\multicolumn{ 1}{|c|}{} &    LHR &   8.1e-4 &         80 &       4000 &        12.7 &   1.3e-3 &         20 &      16031 &        14.1 \\

\hline \hline
       400 &   PCP &   4.3e-1 &         207 &      83954 &        26.1 &   7.0e-1&        214 &      89370 &        35.7 \\

           &    LHR &   8.2e-4 &         160 &      16218 &       63.4 &   1.7e-3 &         40 &      47999 &       54.3 \\

\hline \hline
       800 &   PCP &   4.8e-1 &         414 &      336188 &        36.2 &   9.3e-2 &        348 &     355878 &        48.2 \\

           &    LHR &   9.9e-4 &         20 &      64283 &       91.7 &   2.1e-3 &         320 &     191998 &       108.2 \\
\hline
\end{tabular}
\label{tab:evaluation}
\end{table*}

From the results, obviously, compared with PCP, LHR model could exactly  recover the matrix from higher ranks and denser errors. However, the table just provides two discrete tests. We will provide more thorough investigation in the next subsection.
\subsubsection{ Feasible region}
\begin{figure*}[!htp]
  \centering
    \subfigure[Feasible region verification.]{
    \label{fig:Tolerance:feasible}
\includegraphics[width=0.45\linewidth]{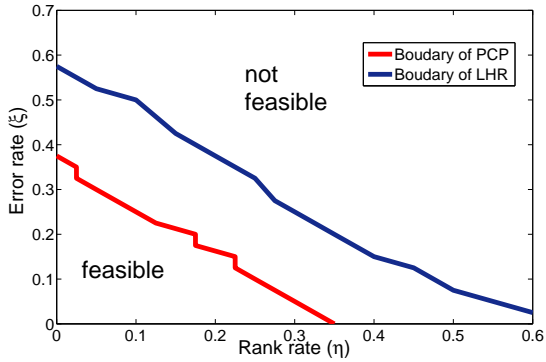}}
    \subfigure[Convergence verification.]{
      \includegraphics[width=0.45\linewidth]{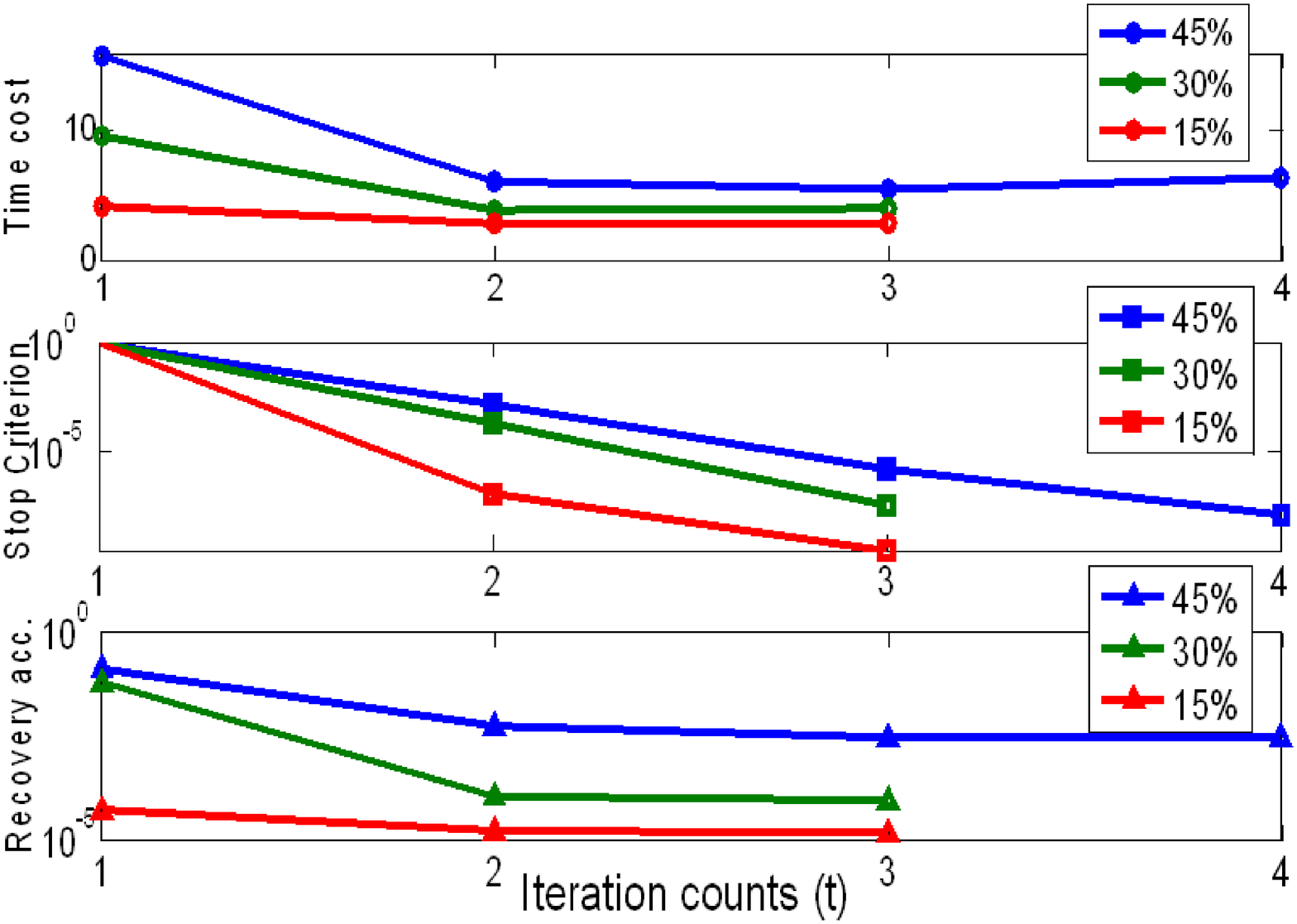}}
    \label{fig:Tolerance:covergence}
  \caption{Feasible region   and the convergence verifications.}
\label{fig:Tolerance}
\end{figure*}
Since the basic optimization involves two terms,
i.e., low rank matrix and sparse error. In this part, we will varies these two variables to test the feasible boundary of PCP and LHR, respectively.
The experiments are conducted on the $400
\times 400 $ matrices with sparse errors uniformly distributed in
$[-100,100]$. In the feasible region verification, when the recovery
accuracy is larger than $1\%$ (i.e.,
$\frac{\|A-A^*\|_F}{\|A^*\|_F}>0.01$), it is believed that the
algorithm diverges. The two rates $\eta$ and $\xi$ are varied from zero to
one with the step of $0.025$. On each test point, both the PCP
and LHR are repeated for $10$ times. If the median recovery accuracy
is less than $1\%$, the point is regarded as the feasible point. The
feasible regions of these two algorithms are shown in Fig.\ref{fig:Tolerance:feasible}.

From Fig.\ref{fig:Tolerance:feasible}, the feasible region of LHR is much
larger than the region of PCP. We  get the same conclusion as
made in \cite{Candes:1229741} that the feasible
boundary of PCP roughly fits the curve that
$\eta^{PCP}+\xi^{PCP}=0.35$.  The boundary of LHR is around the
curve that $\eta^{LHR}+\xi^{LHR}=0.575$. Moreover, on the two sides
of the red curve in Fig.\ref{fig:Tolerance:feasible}, the boundary equation
can be even extended to $\eta^{LHR}+\rho^{LHR}=0.6$. From this
test, it is apparent that the proposed LHR algorithm covers a larger
area of the feasible region, which implies that LHR could handle more
difficult tasks that robust PCA fails to do.

\subsubsection{Convergence discussions}

Finally, we will experimentally verify the convergence of the LHR. The experiments are conducted on $400
\times 400$ matrices with the rank equivalent to $40$ and the
portion of gross errors are set as $15\%,30\%$ and $45\%$,
respectively. The experimental results are reported in Fig.\ref{fig:Tolerance:covergence} where the axis's coordinate denotes the
iteration sequences, i.e. the count $t$ in
Algorithm.\ref{alg:MM-LHR}.


The top sub-figure in Fig.\ref{fig:Tolerance:covergence} reports the time cost of
each iteration. It is interesting to note that the denser the error
is, the more time cost is required for one iteration. Besides,
the most time consuming part occurs in the first iteration. During
the first iteration, (\ref{eqs:rew}) subjects to the
typical PCP problem. However, in the second and the third iteration, the weight matrix is assigned with different
values and thus it could make (\ref{eqs:rew}) converge with less
iterations. Therefore, the time cost for each iteration is different in LHR. The first iteration needs many computational resources while the later
ones can be further accelerated owing to the penalty of the
weight matrix.

The middle sub-figure records the stopping criterion, which is
denoted as $\frac{\|W^{(t+1)}-W^{(t)}\|_F}{\|W^{(t)}\|_F}$. It is
believed that the LHR converges when the stopping criterion is less
than $1e-5$. It is apparent from Fig.\ref{fig:Tolerance:covergence} that the LHR
could converge in just three iterations with $15\%$ and $30\%$ gross
errors. While for the complicated case with $45\%$ errors, LHR can
converge in four steps. The bottom sub-figure shows the recovery accuracy after each
iteration. It is obvious that the recovery accuracy increases
significantly from the first iteration to the second one.
Such an increase phenomenon verify the advantage of the
reweighted approach derived from LHR.

\subsection{Practical applications}
PCP is a powerful tool for many practical applications. Here, we will conduct two practical applications to verify the effectiveness of PCP and LHR on real-world data.
\subsubsection{Shadow and specularities removal from faces}
Following the framework suggested in \cite{Candes:1229741}, we stack the faces of the same subject under different lighting conditions  as the columns in a matrix $\mathbf{P}$. The experiments are conducted on extended Yale-B dataset where each face is with the resolutions of $192 \times 168$. Then, the corrupted matrix $\mathbf{P}$ is recovered by PCP and LHR, respectively. After recovery, the shadows, specularities and other reflectance are removed in the error matrix ($\mathbf{E}$)  and the clean faces are accumulated in the low rank matrix ($\mathbf{A}$).
\begin{figure}[!htp]
  \centering
    \subfigure[Dense shadow]{
    \label{fig:faceshadow:DenseShadow}
    \includegraphics[width=0.8\linewidth]{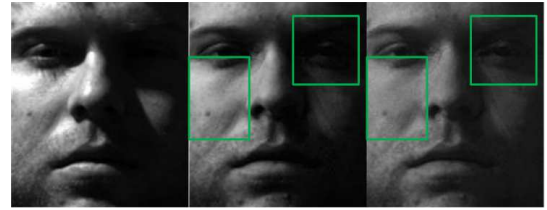}}
    \subfigure[Shadow texture]{
    \label{fig:faceshadow:texture}
    \includegraphics[width=0.8\linewidth]{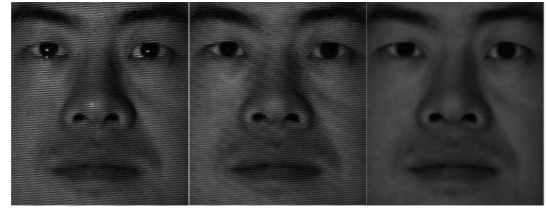}}
  \caption{Shadow and specularities  removal from faces (best viewed on screen).}
\label{fig:faceshadow}
\end{figure}
The experimental results are provided in Fig.\ref{fig:faceshadow}, where in each sub-figure from left to right are: original faces in Yale-B (left), faces recovered by PCP (median) and faces recovered by LHR (right), respectively. It is greatly recommended to enlarge the faces in Fig.\ref{fig:faceshadow} to view the details. In Fig.\ref{fig:faceshadow:DenseShadow}, when there exist dense shadows on the face image, the effectiveness of LHR becomes apparent to remove the dense shadows distribute on the left face. The dense texture removal ability is especially highlighted in Fig.\ref{fig:faceshadow:texture}, where there are significant visual contrasts  between the faces recovered by PCP and LHR. The face recovered by LHR is much clean.
\subsubsection{Video surveillance}
\begin{figure}[!htp]
  \centering
  \subfigure[HW(439 frames).]{
    \includegraphics[width=0.8\linewidth]{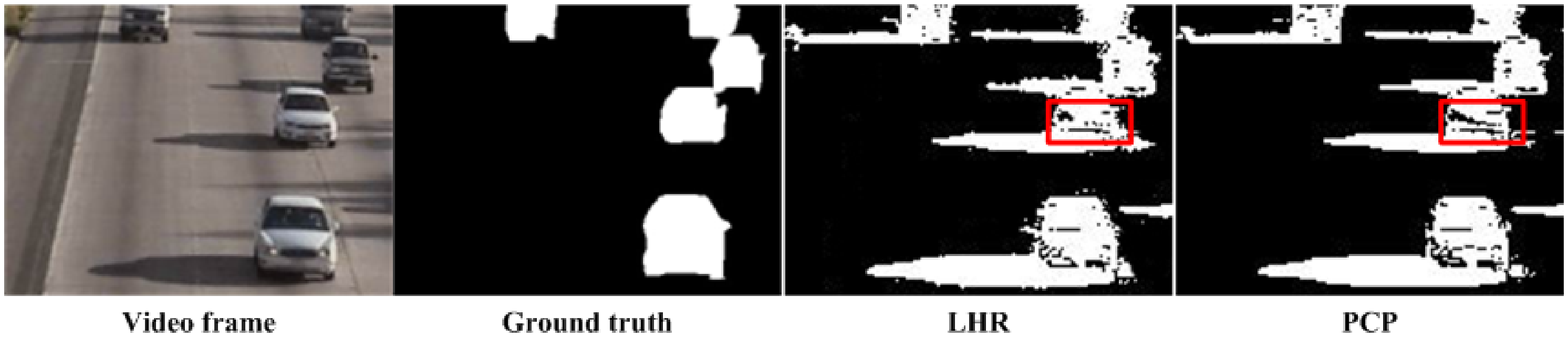}}
    \subfigure[Lab(886 frames)]{
    \includegraphics[width=0.8\linewidth]{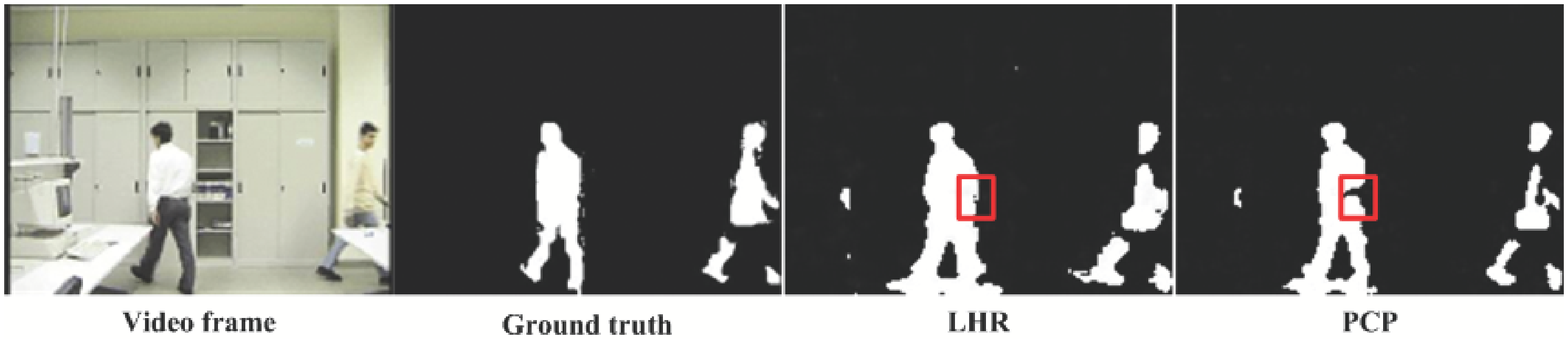}}
    \subfigure[Seam(459 frames)]{
   \label{fig:seam}
    \includegraphics[width=0.8\linewidth]{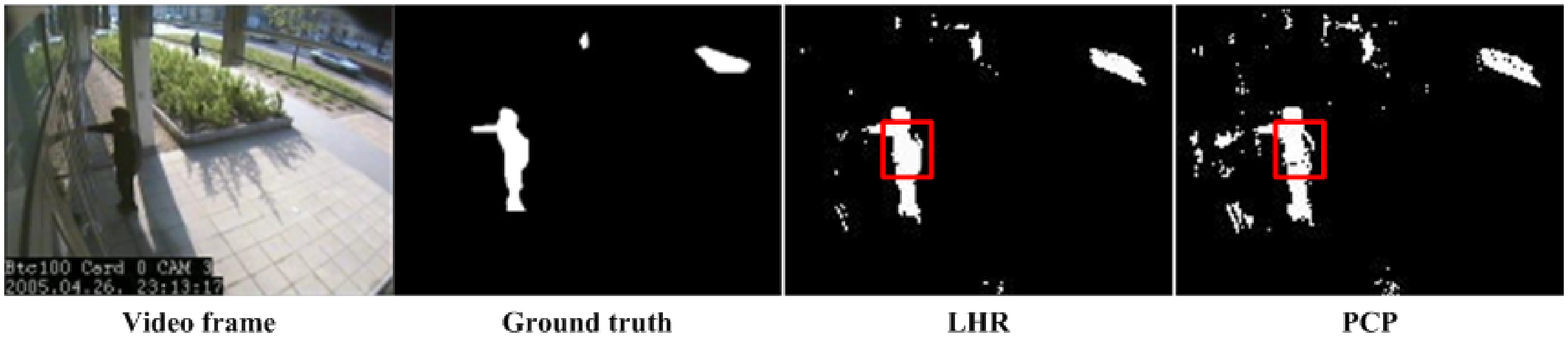}}
  \caption{Benchmark videos for background modeling. In each sub-figure,  from left to right are original video frames, foreground ground truth, LHR result and PCP result,respectively.}
\label{fig:VS}
\end{figure}

The background modeling can also be categorized as a low rank matrix recovery problem, where the backgrounds correspond to the low rank matrix $\mathbf{A}$ and the foregrounds are removed in the error matrix $\mathbf{E}$. We  use the  videos and ground truth in \cite{VS} for quantitative evaluations. Three videos used in this experiment are listed in Fig.\ref{fig:VS}.

For the sake of computational efficiency, we normalize each image to the resolutions of $120 \times 160$ and all the frames are converted to gray-scaler. The benchmark videos used here contain too many frames which lead to a large matrix. It is theoretical feasible to use the two methods for any large matrix recovery. Unfortunately, for practical implementation, large matrices are always beyond the memory limitation of Matlab. Therefore, for each video, we uniformly divide the large matrix to be sub-matrices which has less than 200 columns. We recover these sub-matrices by setting $\lambda=\frac{1}{\sqrt{m}}$, respectively.

The segmented foregrounds and the ground truth are shown in Fig.\ref{fig:VS}. From the results we know that LHR could remove much denser errors from the corrupted matrix rather than PCP. Such claim is verified from three sequences in Fig.\ref{fig:VS} that LHR makes much complete object recovery from the video. Besides, in Fig.\ref{fig:seam}, it is also apparent that LHR only keeps dense errors in the sparse error term. In the seam sequences, there are obvious illumination changes in different frames. PCP is sensitive to these small variations and thus makes much more small isolated noise parts in the foreground. On the other hand, LHR is much robust to these local variations and only keeps dense corruptions in the sparse term.

Although there are many advanced techniques for video background modeling, it is not the main concern of this work. Therefore, without the loss of generality, we use the Mixture of Gaussian (MoG)\cite{MOG_VS} as the comparison baseline. For evaluation, both the  false negative rate (FNR) and false positive rate (FPR) are calculated in the sense of foreground detection. These two scores exactly correspond to the Type I and Type II errors in machine learning, whose definitions may refer to \cite{FNR}.  FNR indicates the ability of the method to correctly recover the foreground and FPR represents the potential of a method on distinguishing the background. Both these two rates are judged by the criterion that the less the better. The experimental results are tabulated in table.\ref{tab:vsresult}. We also report the time cost (in minutes) of PCP and LHR on these videos. But we omit the time cost of MoG since it can be finished in almost real time.

\begin{table}
\centering
\label{tab:vsresult}
\caption{Quantitative evaluation of PCP and LHR for video surveillance. }
\begin{tabular}{|p{0.45cm}|p{0.45cm}p{0.45cm}p{0.45cm}|p{0.45cm}p{0.45cm}p{0.45cm}|p{0.48cm}p{0.48cm}|}
\hline
\multicolumn{ 1}{|c}{Data} & \multicolumn{3}{|c}{False Negative Rate\%} & \multicolumn{ 3}{|c}{False Positive Rate\% } & \multicolumn{ 2}{|c|}{Time($m$)} \\

\multicolumn{ 1}{|c|}{} &        MoG &       PCP &        LHR &        MoG &       PCP&        LHR &       PCP &        LHR \\
\hline\hline
  HW &       22.2 &       18.7 &       14.3 &        8.8 &        7.8 &        8.4 &       13.2 &       23.5 \\
\hline
Lab. &       15.1 &       10.1 &        8.3 &        6.7 &        6.4 &        6.1 &       25.4 &       43.7 \\
\hline
      Seam &       23.5 &       11.3 &        9.2 &        9.7 &        6.1 &        6.3 &       11.4 &       19.9 \\
\hline
\end{tabular}

\end{table}

From the results, PCP and LHR greatly outperform the performance of MoG. Moreover, LHR has lower FNRs than PCP which implies that LHR could better detect the foreground than PCP.  However, on the video highway and seam, the FPR score of LHR is a little worse than PCP. One possible reason may ascribe to that there are too many moving shadows in these two videos, where both the objects and shadows are regarded as errors. In the ground truth frames, the shadows are regarded as background. LHR  could recover much denser errors form a low rank matrix  and thus causes a comparable low FNR score.

\section{LHR for low rank representation}
\label{sec:REWLRR}
In this part, LHR will be applied to the task of low rank representation (LRR)\cite{LRR}\cite{ICML-LRR} by formulating the constraint as $\mathbf{P}=\mathbf{PA}+\mathbf{E}$, where the correlation affine matrix $\mathbf{A}$ is low rank and the noises in $\mathbf{E}$ are sparse. In the remaining parts of this section, we will first show how to use the joint optimization strategy to solve the LRR problem by LHR model. Then, two practical applications on motion segmentation and stock clustering will be presented and discussed.
\subsection{Joint optimization strategy of LHR for LRR}
When applying LHR to low rank representation, we should solve a sequence of convex optimizations in the form,
\begin{equation}
\begin{array}{l}
 \min.~\left\| {{\bf{W}}_Y {\bf{AW}}_Z } \right\|_*  + \lambda \left\| {{\bf{W}}_E  \odot {\bf{E}}} \right\|_{\ell _1 }  \\
 s.t.~~~~{\bf{P}} = {\bf{PA}} + {\bf{E}} \\
 \end{array}
 \end{equation}
 To make the nuclear norm trackable, we add an equality and tries to solve
\begin{equation}
\begin{array}{l}
 \min .\left\| {\bf{J}} \right\|_*  + \lambda \left\| {{\bf{W}}_E  \odot {\bf{E}}} \right\|_{\ell _1 }  \\
 s.t.~~~{\bf{b}}_1  = {\bf{P}} - {\bf{PA}} - {\bf{E}} = {\bf{0}} \\
 ~~~~~~~{\bf{b}}_2  = {\bf{J}} - {\bf{W}}_Y {\bf{AW}}_Z  = {\bf{0}} \\
 \label{eqs:rewlrr}
 \end{array}
 \end{equation}
Using the ADM strategy and following the similar derivations introduced in subsection \ref{sec:RPCA-solve}, we can solve the optimization in (\ref{eqs:rewlrr}) and we  directly provide the update rules for each variable in algorithm \ref{alg:rewlrr}.
\begin{algorithm}
\label{alg:rewlrr}
$E_{ij}^{k}  = s_{\lambda \mu ^{ - 1} \left| {(W_E^{(t)})_{ij} } \right|} (P - PA^{k-1}
- \mu ^{ - 1} C_1^k)_{ij},\forall ij$\;
          $\mathbf{J}^{k}  = d_{\mu^{-1}}  ({\bf{W}}_{Y}^{(t)} {\bf{A}^{k-1}\mathbf{W}}^{(t)}_{Z}  + \mu ^{ - 1} {\bf{C}}_2^k )$\;
          ${\bf{A}^k}  = {\bf{A}}^{k-1}  + \gamma [{\bf{W}}_Y^{(t)} ({\bf{b}}_1^k  + \mu ^{ - 1} {\bf{C}}_2^k ){\bf{W}}_Z^{(t)} + \mathbf{P}^T({\bf{b}}_2^k  + \mu ^{ - 1} {\bf{C}}_1^k )]$\;
          \tcp{Dual ascending.}
          $\mathbf{C}_1^k=\mathbf{C}_1^{k-1}+\mu_k \mathbf{b}_1^k$\;
          $\mathbf{C}_2^k=\mathbf{C}_2^{k-1}+\mu_k \mathbf{b}_2^k$\;
\caption{Update rule for  the variables in (\ref{eqs:rewlrr})}
\label{alg:ADM}
\end{algorithm}
To show LHR ideally represents low rank structures from data, experiments on subspace clustering are conducted on two datasets. First, we test LHR on slightly corrupted data-set, i.e., Hopkins156
motion database. Since the effectiveness of LHR are especially emphasized on the data with great corruptions. We will also consider one practical application of using LHR for stock clustering.
\subsection{Motion segmentation in video sequences}

In this part, we apply LHR to the task of motion segmentation in Hopkins155 dataset \cite{SC}. Hopkins155 database is a benchmark platform to evaluate general subspace clustering algorithms, which contains 156 video sequences and each of them has been summarized to be a matrix recoding 39 to 50 data vectors. The primary task of subspace clustering is to categorize each motion to its corresponding subspace, where  each video corresponds to	
a sole clustering task and it leads to 156 clustering tasks in total.

For comparisons, we will compare LHR with LRR as well as other benchmark algorithms for subspace clustering. The comparisons include Random Sample Consensus (RANSAC)\cite{RANSAC}, Generalized Principal Component Analysis (GPCA)\cite{GPCA}, Local Subspace Affinity (LSA), Locally Linear Manifold Clustering (LLMC) and Sparse Subspace Clustering (SSC). RANSAC is a statistic method which clusters data by iteratively distinguishing the data by inliers and outliers. GPCA presents an algebraic method to cluster the mixed data by the normal vectors of the data points. Manifold based algorithms, e.g. LSA and LLMC, assume that one point and its neighbors span as a linear subspace and they are clustered via spectral embedding. SSC assumes that the affine matrix between data are sparse and it segments the data via normalized cut\cite{normalizedcut}.

 In LRR \cite{LRR}, Liu \emph{et al. } introduced two models that respectively used $\ell_1$ norm  and $\ell_{2,1}$ norm to  penalize sparse corruptions. In this paper, we will only report the results with the comparison to $\ell_{2,1}$ norm since it always performs better than $\ell_1$ penalty in LRR. In order to provide a thorough comparison with LRR, we strictly follow the steps and the default parameter settings  suggested in \cite{LRR}. For LHR model, we choose parameter $\lambda=0.4$. In the experiments of LRR for motion segmentation, some post-processing  are performed on the learned low rank structure to seek for the best clustering accuracy. For example, in LRR, after  getting the representation matrix $\mathbf{A}$, an extra PCP processing are implemented on $\mathbf{A}$ to enhance the low-rankness and such post-processing definitely  increases  SC accuracy. However, the main contribution of this work only focus LHR model on low-rank structure learning while not on the single task of subspace clustering. Therefore, we exclude all the post-processing steps to emphasize the effectiveness of the LRSL model itself. In our result, all the methods are implemented with the same criterion to avoid bias treatments.

Hopkins155 contains two subspace conditions in a video sequence, i.e., with two motions or three motions  and thus we report the segmenting errors  for two subspaces (TWO), three subspaces (THREE) and for both conditions (ALL) in Table.\ref{tab:SGResults}. From the results we know that sparse based methods generally outperform other algorithms for motion segmentation. Among three sparse methods, LHR gains the best clustering accuracy. However, the accuracy only has slight improvements on LRR. As indicated in \cite{LRR}, motion data only contains small corruptions and LRR could already achieve promising performance with the accuracy higher than 90\%. With some post-processing implementations, the accuracy can even be further improved. Therefore, in order to highlight the effectiveness of LHR on low rank representation with corrupted data, some more complicated problems will be considered.

 \begin{table}[!htbp]
   \label{tab:SGResults}
  \centering
  \caption{Motion segmentation errors (mean ) of several
algorithms on the Hopkins155 motion segmentation database.}
\begin{sc}
\begin{tabular}{|c|c|ccc|}
\hline
   Category & Method&  TWO &THREE &ALL\\
\hline
 \emph{Algebraic}  &   GPCA & 11.2  & 27.7  & 14.2 \\
\hline
 \emph{Statistic}  &   RANSAC & 8.9  & 24.1  & 12.5  \\
\hline
\multicolumn{ 1}{|c|}{\emph{Manifold}}&LSA & 8.7  & 21.4  & 11.6  \\
\multicolumn{ 1}{|c|}{}& LLMC  & 8.1 & 20.8  & 10.9  \\
\hline
\multicolumn{ 1}{|c|}{}& SSC  & 5.4 & 15.3  & 7.6 \\
\multicolumn{ 1}{|c|}{\emph{Sparse}} & LRR & 4.4  & 14.9  & 6.7  \\
\multicolumn{ 1}{|c|}{} &        \textbf{LHR} & 3.1  & 13.9  & 5.6 \\
\hline
\end{tabular}
\end{sc}
\end{table}

%

\subsection{Stock clustering}
It is not trivial to consider applying LHR model to more complicated practical data where the effectiveness of LHR  on corrupted data will be over emphasized. In practical world, one of the most difficult data structures to be analyzed is the stock price which can be greatly affected by company news, rumors and global economic atmosphere. Therefore, data mining approaches of financial signals have been proven to be very difficult but on the other hand, it is very profitable.

In this paper, we will discuss how to use the LRR and LHR model to the interesting, albeit not very lucrative, task of stock clustering based on their industrial categories. In many stock exchange centers around the world,  stocks are always divided into different industrial categories. For example, on the New York Stock Exchange Center, IBM and J.P.Morgan are respectively categorized into the\emph{ computer based system}  category and \emph{money center banks} category. It is generally assumed that stocks in the same category always have similar market performance. This basic assumption is widely used by many hedge funds for statistic arbitrage. In this paper, we consider that stocks in the same industrial category span as a subspace and therefore the  goal of stock clustering, a.k.a. stock categorization, is to identify a stock's industrial label by its historical prices.

The experiments are conducted on stocks from two global stock exchange markets in New York and Hong Kong.  In each market, we choose 10 industrial categories which have the largest market capitalizations. The categories divided by the exchange centers  are used as the ground truth label. In each category, we only choose the stocks whose market capitalizations are within the top 10 ranks in one category. The stock prices on New York market are obtained from \cite{yahoofinance} and the stock prices in Hong Kong market are obtained from \cite{googlefinance}. Unfortunately, some historical prices for stocks in \cite{yahoofinance} are not recorded and provided\footnote{For example, in the industrial category of Drug Manufactures, it is not possible to get the historical data of CIPILA.LTD from \cite{yahoofinance} which is the only the interface for us to get the stock prices in US. }. Therefore, for the US market, we  accumulated 76 stocks divided into 10 classes and each class contains 7 to 9 stocks; for Hong Kong market, we obtain 96 stocks spanning 10 classes. For classification, the weekly closed prices  from $07/01/2008$ to $31/10/2011$ including $200$ weeks,  are used because financial experts always look at weekly close prices to judge the long-term trend of a price.

\begin{figure}[!htp]
  \centering
    \includegraphics[width=1\linewidth]{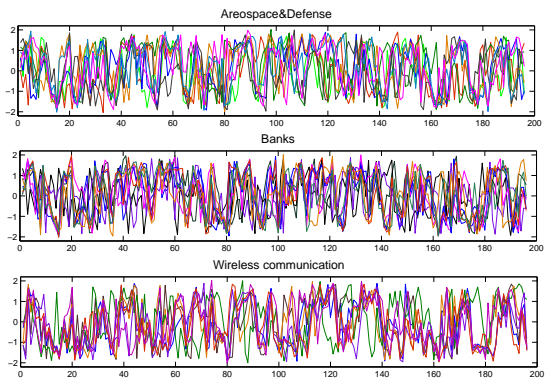}
  \caption{Normalized stock prices in NY of the categories: Areospace\&Defense,Banks and Wireless communication. In each category, lines in different colors represent different stocks. (best viewed on screen)}
\label{fig:stock}
\end{figure}

As stated previously, the stock prices may have sharp drop and up which cause outliers in the raw data. Besides, the prices of different stocks are various that cannot be evaluated with the same quantity scaler. For the ease of pattern mining, we use the time-based normalization strategy suggested in \cite{stockmining}\cite{financetimeseries} to pre-process the stock prices:
\[\tilde p(t) = \frac{{p(t) - \mu _\alpha  (t)}}{{\sigma _\alpha  (t)}},
\]
where $p(t)$ is the price of a certain stock at time $t$ , $\mu_\alpha(t)$ and
$\sigma _\alpha(t)$ are respectively the average value and standard derivation of the stock
prices in the interval $[t-\alpha,t]$. We plot the normalized stock prices of three categories in Fig.\ref{fig:stock}. After normalization, we further adopt PCA method to reduce the dimensions of stocks from $\mathbb{R}^{200}$ to $\mathbb{R}^5$. Theoretically, the rank of subspaces after PCA should be $10-1=9$ because it contains 10 subspaces and the rank is degraded by 1 during PCA implementation. But, in the simulation, we find that the maximal clustering accuracies for both markets are achieved with the PCA dimensions of 5.

\begin{table}
\centering
  \caption{Clustering errors of the stocks in ten categories from New York and Hong Kong markets.}
    \begin{tabular}{ccccc}
    \toprule
    Markets & GPCA  & RANSAC & LSA   & LLMC \\
    \midrule
    New York & 60.1  & 59.3  & 51.7  & 54.3 \\
    Hong Kong & 57.3 & 55.8  & 54.7   &53.7  \\
    \hline
    \hline
    Markets & SSC   & LRR   & LHR   &  \\
    New York & 48.6  & 44.1  & 36.2  &  \\
    Hong Kong & 49.1  & 46.5  & 38.3  &  \\
    \bottomrule
    \end{tabular}%
  \label{tab:stockresult}%
\end{table}%

The clustering errors of different SC methods on the stocks from these two markets are summarized in Table.\ref{tab:stockresult}. From the results, it is obvious that LHR significantly outperforms other methods. It improves statistic and graph based methods for about $20\%$. Among all the sparse methods, LHR makes improvements on LRR for about $8\%$. Although LHR performs the best among all the methods, the clustering accuracy is only about $63\%$ and $61\%$ on US and Hong Kong markets, respectively. The clustering accuracy is not as high as those on the motion data. This may be ascribe to that the raw data  and ground truth label themselves contain many uncertainties. See the bottom sub-figure in Fig.\ref{fig:stock} for the stocks in the \emph{wireless communication} category, the normalized stock marked with the green color performs quite different from other stocks in the same category. But the experimental results reported here is sufficient to verify the effectiveness of subspace clustering for 10 classes categorization. If no intelligent learning approaches were imposed, the expected accuracy may be only $10\%$. Although with such "bad" raw data, the proposed LHR could achieve the accuracy as high as $62\%$ in a definitely unsupervised way.
\section{Conclusion}
\label{sec:Con}
This paper presents a log-sum heuristic recovery algorithm to learn the essential low rank structures from corrupted matrices. We introduced a MM algorithm to convert the non-convex objective function a series of convex optimizations via reweighed approaches and proved that the solution may converge to a stationary point. Then, the general model was applied to two practical tasks of LRMR and SC. In both of the two models, we gave the solution/update rules to each variable in the joint optimizations via ADM. For the general PCP problem, LHR extended the feasible region to the boundary of $\eta+\xi=0.58$. For SC problem, LHR achieved state-of-the-art results on motion segmentation and achieved promising results on stock clustering which contain too many outliers and uncertainties.  However, a limitation of the proposed LHR model is for the reweighted phenomenon that requires to solve  convex optimizations for multiple times. The implementations of LHR is a bit more time consuming than PCP and LRR. Therefore, LHR model is especially recommend to learn the low rank structure from data with denser corruptions and higher ranks.

\appendix
\subsection{ABBREVIATIONS}
\textbf{ADM}: Alternating Direction Method, \textbf{MM}: Majorization Minimization,  \textbf{GPCA}:Generalized Principal Component Analysis, \textbf{PCP}: Principal Component Pursuit, \textbf{LHR}:Log-sum Heuristic Recovery, \textbf{LLMC}: Locally Linear Manifold Clustering, \textbf{LRMR}: Low Rank Matrix Recovery, \textbf{LRR}:Low rank representation, \textbf{LRSL}: Low Rank Structure Learning, \textbf{LSA}: Local Subspace Affinity,\textbf{MoG}: Mixture of Gaussian, \textbf{RANSAC}: Random Sample Consensus, \textbf{RPCA}:Robust Principal Component analysis,  \textbf{SC}:Subspace Clustering, \textbf{SSC}: Sparse Subspace Clustering.

\subsection{Convergence of subsequences in the proof of  lemma \ref{Lemma:Convergence}}
\label{sec:subconvergent}
In this part, we provide the discussions about the properties of the convergent subsequences that are used in the proof of Lemma \ref{Lemma:Convergence}.

Since sequence $\hat{X^{t}}=\{\mathbf{Y}^{t},\mathbf{Z}^{t},\mathbf{A}^{t},\mathbf{E}^{t}\}$ is generated via  Eq.\ref{eqs:LHR}, we know that $\hat{X}\in\hat{D}$  strictly holds. Therefore, all the variables (i.e. $\mathbf{Y}^{t},\mathbf{Z}^{t},\mathbf{A}^{t},\mathbf{E}^{t}$) in set $\hat{X}$ should be bounded. This claim can be easily verified because that if any variable in the set $\hat{X}$ goes to infinity, the constraints in domain $\hat{D}$ will not be satisfied. Accordingly, we know that sequence $\hat{X}^{t}$ is bounded.
 According to the \emph{Bolzano-Welestrass Theorem} \cite{BW-theorem}, we know that  every bounded sequence has a convergent subsequence. Since $\hat{X}^{t}$ is bounded, it is apparent that there exists a convergent subsequence $\hat{X}^{t_k}$. Without the loss of generality, we can construct another subsequence $\hat{X}^{t_k+1}$ which is also convergent.  The proof of the convergence of $\hat{X}^{t_k+1}$ relies on the monotonically decreasing property proved in Lemma \ref{Lemma:Decrease}. Since $H(\cdot)$ is monotonically decreasing, it is easy to check that$
H(\hat{X}^{t_k } ) \ge H(\hat{X}^{t_k  + 1} ) \ge H(\hat{X}^{t_{k + 1} } ) \ge H(\hat{X}^{t_{k + 1}  + 1} ) \ge H(\hat{X}^{t_{k + 1 + 1} } )$.
According to the above inequalities, we get that,
\begin{equation}
\mathop {\lim }\limits_{k \to \infty } H(\hat{X}^{t_k } ) \ge \mathop {\lim }\limits_{k \to \infty } H(\hat{X}^{t_{k + 1} } ) \ge \mathop {\lim }\limits_{k \to \infty } H(\hat{X}^{t_{k + 2} } )
\label{eqs:Seq}
\end{equation}
Since subsequence $\hat{X^{t_k}}$ converges, it is obvious that $
\mathop {\lim }\limits_{k \to \infty } H(\hat{X}^{t_k } )=\mathop {\lim }\limits_{k \to \infty } H(\hat{X}^{t_{k+2} } )=\beta $. According to the famous \emph{Squeeze Theorem }\cite{Sequezz-theorem}, from (\ref{eqs:Seq}), we get the $\mathop {\lim }\limits_{k \to \infty } H(\hat{X}^{t_{k}+1 } )=H(\mathop {\lim }\limits_{k \to \infty }\hat{X}^{t_{k}+1 })=\beta $. Since the function $H(\cdot)$ is monotonically decreasing and $\hat{X}$ is bounded, the convergence of  $H(\hat{X}^{t_{k}+1 } )$ can be obtained if and only if the subsequence $\hat{X}^{t_{k}+1}$ is convergent.

%
%

\ifCLASSOPTIONcaptionsoff
  \newpage
\fi

\bibliographystyle{ieeetran}
\bibliography{instructions2}

%

%
%
%




\end{document}